\documentclass[11pt]{amsart}
\usepackage{amssymb}
\usepackage{amsmath}

\newcommand{\mpn}{\medskip\par\noindent}
\newcommand{\pn}{\par\noindent}

\newcommand{\bpn}{\bigskip\par\noindent}
\theoremstyle{definition}

\theoremstyle{remark}

\numberwithin{equation}{section}

\begin{document}
\newcommand{\Mod}[1]{\,(\text{\mbox{\rm mod}}\;#1)}
\title[On the associated sequences of special polynomials ]
{On the associated sequences of special polynomials }
\author{D. S. Kim, T. Kim, and S.-H. Rim}
\begin{abstract}
 In this paper, we investigate some properties of the associated sequence of Daehee and Changhee polynomials. Finally, we give some interesting identities of associated sequence involving some special polynomials.
\end{abstract}

 \maketitle

\section{\bf Introduction}
In this paper, we assume that
 $\lambda \in \mathbb{C}$, with $\lambda \neq 1 $. For $ \alpha\in \mathbb{R}$, the Frobenius
-Euler polynomials of order $r$  are defined by the generating function to be
\begin{align*}\tag{1}
\left(\frac{1-\lambda}{e^t-\lambda}\right)^{\alpha}e^{xt}=\sum_{n=0}^{\infty}H_n^{(\alpha)}(x
| \lambda)\frac{t^n}{n!}, \quad ({see \,[9,11,13,22,27]}).
\end{align*}
In the special case, $x=0$, $H_n^{(\alpha)}(0| \lambda)=H_n^{(\alpha)}(\lambda)$ are called $n$-th Frobenius-Euler numbers of order $\alpha$, $({see \,[8,9,27]})$.
As is well known, the Bernoulli polynomials of order $\alpha$ are defined by the generating function to be:
\begin{align*}\tag{2}
\left(\frac{t}{e^t-1}\right)^{\alpha}e^{xt}=\sum_{n=0}^{\infty}B_n^{(\alpha)}(x)\frac{t^n}{n!}, \quad ({see \,[10,12,21,28]}).
\end{align*}
In the special case, $x=0$, $B_n^{(\alpha)}(0| \lambda)=B_n^{(\alpha)}$ are called $n$-th Bernoulli numbers of order $\alpha$.
The Stirling numbers of the second kind are defined by
\begin{align*}\tag{3}
(e^t-1)^n=n!\sum_{l=n}^{\infty}S_2(l,n)\frac{t^l}{l!},\quad ({see \,[5,20,21,26]}).
\end{align*}
and the Stirling numbers of the first kind are given by
\begin{align*}\tag{4}
(x)_n=x(x-1)...(x-n+1)=\sum_{l=0}^{n}S_1(l,n)x^l, \quad({see \,[5,20,21]}).
\end{align*}
Let $\mathcal{F}$ be the set of all formal power series in the variable $t$ over $\mathbb{C}$ with
\begin{align*}\tag{5}
\mathcal{F}=\{f(t)=\sum_{n=0}^{\infty}\frac{a_k}{k!}t^k ~|~a_k \in\mathbb{C}\}.
\end{align*}
Let us assume that $\mathbb{P}$ is the algebra of polynomials in the variable $x$ over $\mathbb{C}$  and  $\mathbb{P^*}$ is the vector
space of all linear functionals on $\mathbb{P}$.
$<L ~|~ p(x)>$ denotes the action of the linear functional
$L$ on a polynomial $p(x)$ and we remind that the vector space structure on $\mathbb{P^*}$ is defined by
$<L+M ~|~p(x)>=\\<L~|~p(x)>+<M~|~p(x)>,$
 and  $<cL ~|~ p(x)>=c<L~|~p(x)>,$
\\where $c$ is a complex constant, $(see [6,9,20,21])$.
For $f(t)=\sum_{k=0}^{\infty}\frac{a_k}{k!}t^k \in \mathcal{F}$, we define a linear functional on $\mathbb{P}$ by setting
\begin{align*}\tag{6}
<f(t) ~ | ~x^n>=a_n ,(n\geq 0),\quad (see [20,21]).
\end{align*}
From (6), we note that
\begin{align*}\tag{7}
<t^k ~ | ~x^n>=n!\delta_{n,k}, \quad (n ,k \geq 0),
\end{align*}
where $\delta_{n,k}$ is the Kronecker's symbol.
Let $f_L(t)=\sum_{k=0}^{\infty}\frac{<L~|~ x^k>}{k!}t^k $. Then, by (7), we see that
$<f_L(t) ~|~ x^n> = <L~|~x^n>$.
The map $L \mapsto f_L(t)$ is a vector space isomorphism
from $\mathbb{P^*}$ onto  $\mathcal{F}$. Henceforth, $\mathcal{F}$
is thought of as both a formal power series and a linear functional. We call $\mathcal{F}$
the umbral algebra. The umbral calculus is the study
of umbral algebra $(see [7,11,,19,20,21])$.
The order $O(f(t))$ of the non-zero power series $f(t)$ is the smallest integer $k$ for which the coefficient of $t^k$ does not vanish $(see [10,20,21])$. If $O(f(t))=1$, then $f(t)$
is called a delta series  and a series $f(t)$ having $O(f(t))=0$ is called an invertible series $(see[10,20,21])$. Let $f(t)$ be a delta series and $g(t)$ be an invertible series. Then there exists a unique sequence $S_n(x)$ of polynomials such that $<g(t)f(t)^k |S_n(x)>=n!\delta_{n,k}$ where $n,k \geq 0$. The sequence $S_n(x)$ is called the Sheffer sequence for $(g(t),f(t))$ which is denoted by $S_n(x) \sim (g(t),f(t))$. If $S_n(x)\sim(1,f(t))$, then $S_n(x)$ is called the associated sequence for $f(t)$. By (7), we see that $<e^{yt}~ |~ p(x)>=p(y) $. For $f(t) \in \mathcal{F}$ and $p(x) \in \mathbb{P}$, we have
\begin{align*}\tag{8}
f(t)=\sum_{k=0}^{\infty}\frac{<f(t)~|~x^k>}{k!}t^k ,\quad
p(x)=\sum_{k=0}^{\infty}\frac{<t^k~|~p(x)>}{k!}x^k,\quad (see[11,20,21]) .
\end{align*}
and
\begin{align*}\tag{9}
&<f_1(t)f_2(t) ... f_m(t)~|~x^n>
\\&=\sum_{i_1+...+i_m=n}^{}\binom{n}{i_1,...i_m} <f_1(t)~|~x^{i_1}>... <f_m(t)~|~x^{i_m}>,
\end{align*}
where $ f_1(t),f_2(t),  ...,  f_m(t) \in \mathcal{F}$.
\\By (8), we get
\begin{align*}\tag{10}
p^{(k)}(0)= <t^k~|~p(x)> ,\quad <1 ~|~ p^{(k)}(x)>=p^{(k)}(0).
\end{align*}
Thus, by (10), we see that
\begin{align*}\tag{11}
t^kp(x)=p^{(k)}(x)= \frac{d^kp(x)}{dx^k},\quad (k \geq 0), \quad(see[10,20,21]).
\end{align*}
For $S_n(x) \sim (g(t),f(t))$,  we have
\begin{align*}\tag{12}
\frac{1}{g(\bar{f}(t))}e^{y\bar{f}(t)}=\sum_{k=0}^{\infty}\frac{S_k(y)}{k!}t^k, \quad for \quad all \quad y\in \mathbb{C} ,
\end{align*}
where $\bar{f}(t)$ is the compositional inverse of $f(t)$,
\begin{align*}\tag{13}
S_n(x+y)=\sum_{k=0}^{n} \binom {n}{k} p_k(y)S_{n-k}(x)= \sum_{k=0}^{n}\binom {n}{k}S_{n-k}(y) p_k(x),
\end{align*}
where $p_k(y)=g(t)S_k(y) \sim (1,f(t))$, $(see [5,11,20,21])$.
\\Let $p_n(x) \sim (1,f(t))$ and $q_n(x) \sim (1,g(t))$. Then the transfer formula for the associated sequence is implies that
\begin{align*}\tag{14}
q_n(x)=x\left(\frac{f(t)}{g(t)} \right)^n x^{-1}p_n(x),\quad (see[5,11,20,21]).
\end{align*}

In this paper, we investigate some properties of the associated sequence of Daehee and Changhee polynomials arising from umbral calculus. Finally, we derive some interesting identities of associated sequence of special polynomials from (14).

\section{\bf On the associated sequence of special polynomials  }
For $\lambda \in \mathbb{C}$ with $\lambda \neq 1 $, the Daehee polynomials are defined by the generating function to be
\begin{align*}\tag{15}
\sum_{k=0}^{\infty}\frac{D_k(x|\lambda)}{k!}t^k = \left(\frac{(1-\lambda)+t(1+\lambda)}{(1-\lambda)(1-t)} \right)
\left( \frac{1+t}{1-t}\right)^x,
\end{align*}
$(see~[1,2,14,16,17,18])$.
\\From (15), we note that $D_n(x|\lambda)$ is the Sheffer sequence for $ \left(\frac{1-\lambda}{e^t-\lambda},\frac{e^t-1}{e^t+1}\right)$.
That is,
\begin{align*}\tag{16}
D_n(x|\lambda) \sim \left(\frac{1-\lambda}{e^t-\lambda},\frac{e^t-1}{e^t+1}\right).
\end{align*}
\vskip .1in
As is known, the Mittag-Leffler sequence is given by
\begin{align*}\tag{17}
M_n(x) = \sum_{k=0}^{n} \binom {n}{k}(n-1)_{n-k}2^{k}(x)_k  \sim \left(1,\frac{e^t-1}{e^t+1}\right).
\end{align*}
From (16) and (17), we have
\begin{align*}\tag{18}
M_n(x)= \frac{1-\lambda}{e^t-\lambda}D_n(x|\lambda) \sim \left(1,\frac{e^t-1}{e^t+1}\right).
\end{align*}
Let us consider the following associated sequence:
\begin{align*}\tag{19}
S_n(x|\lambda) \sim \left(1,\frac{1-\lambda}{e^t-\lambda}t\right) \quad and \quad  x^n \sim (1,t).
\end{align*}
By (14) and (19), we get
\begin{align*}\tag{20}
S_n(x|\lambda) &=x\left(\frac{\frac{1-\lambda}{e^t-\lambda}t}{t}\right)^{-n}x^{-1}x^n
\\&=(1-\lambda)^{-n}x(e^t-\lambda)^nx^{n-1}
\\&=(1-\lambda)^{-n}x\sum_{l=0}^{n} \binom {n}{l}(-\lambda)^{n-l}e^{lt}x^{n-1}
\\&=(1-\lambda)^{-n}x\sum_{l=0}^{n} \binom {n}{l}(-\lambda)^{n-l}(x+l)^{n-1}.
\end{align*}
Therefore by (20), we obtain the following
theorem.
\par\bigskip
{\bf Theorem 1  }. For $r \in \mathbb{Z}_+ =\mathbb{N}\cup \{0\}$, let $S_n(x|\lambda) \sim (1,\frac{1-\lambda}{e^t-\lambda}t)$. \\Then we have
\begin{align*}
S_n(x|\lambda)=\frac{x}{(1-\lambda)^n}\sum_{l=0}^{n} \binom {n}{l}(-\lambda)^{n-l}(x+l)^{n-1}.
\end{align*}
\par\bigskip
From (14), (17) and (19), we have
\begin{align*}\tag{21}
&M_n(x)\\&=x\left(\frac{\frac{1-\lambda}{e^t-\lambda}t}{\frac{e^t-1}{e^t+1}}\right)^{n}x^{-1}S_n(x)
=x\left(\frac{1-\lambda}{e^t-\lambda}\right)^n\left(\frac{t(e^t+1)}{e^t-1}\right)^nx^{-1}S_n(x)
\\&=x\left(\frac{1-\lambda}{e^t-\lambda}\right)^n\left(t+\frac{2t}{e^t-1}\right)^n x^{-1}S_n(x)
\\&=x\left(\frac{1-\lambda}{e^t-\lambda}\right)^n\sum_{l=0}^{n} \binom {n}{l}t^{n-l}\left(\frac{2t}{e^t-1}\right)^l x^{-1}S_n(x)
\\&=x\left(\frac{1-\lambda}{e^t-\lambda}\right)^n\sum_{l=0}^{n} \binom {n}{l}t^{n-l}\left(\frac{2t}{e^t-1}\right)^l
\frac{1}{(1-\lambda)^n}\sum_{j=0}^{n} \binom {n}{j}(-\lambda)^{n-j}(x+j)^{n-1}
\\&=x\left(\frac{1-\lambda}{e^t-\lambda}\right)^n\sum_{l=0}^{n}\sum_{j=0}^{n} \binom {n}{l}(1-\lambda)^{-n}\binom {n}{j}
(-\lambda)^{n-j}t^{n-l}\left(\frac{2t}{e^t-1}\right)^l (x+j)^{n-1}.
\end{align*}
By (1), (2) and (12), we get
\begin{align*}\tag{22}
B_n^{(\alpha)}(x) \sim \left(\left(\frac{e^t-1}{t}\right)^{\alpha},t\right)\quad H_n^{(\alpha)}(x |\lambda) \sim \left(\left(\frac{e^t-\lambda}{1-\lambda}\right)^{\alpha},t\right),
\end{align*}
and, by (13), we see that
\begin{align*}\tag{23}
B_n^{(\alpha)}(x)=\sum_{k=0}^{n} \binom {n}{k}B_{n-k}^{(\alpha)}x^k,\quad
 H_n^{(\alpha)}(x |\lambda)=\sum_{l=0}^{n} \binom {n}{l}H_{n-l} ^{(\alpha)}x^l.
\end{align*}
Thus, by (21), (22) and (23), we get
\begin{align*}\tag{24}
& M_n(x)
\\&=x\left(\frac{1-\lambda}{e^t-\lambda}\right)^n\sum_{l=0}^{n}\sum_{j=0}^{n} \binom {n}{l}(1-\lambda)^{-n}\binom {n}{j}
(-\lambda)^{n-j}2^lt^{n-l}B_{n-1}^{(l)}(x+j)
\\&=x\left(\frac{1-\lambda}{e^t-\lambda}\right)^n\sum_{l=1}^{n}\sum_{j=0}^{n} \binom {n}{l}(1-\lambda)^{-n}\binom {n}{j}
(-\lambda)^{n-j}2^l(n-1)_{n-l}B_{l-1}^{(l)}(x+j)
\\&=x\left(\frac{1-\lambda}{e^t-\lambda}\right)^n \sum_{l=1}^{n}\sum_{j=0}^{n} \binom {n}{l}(1-\lambda)^{-n}\binom {n}{j}
(-\lambda)^{n-j}2^l \frac{(n-1)!}{(l-1)!}\sum_{m=0}^{l-1} B_{l-1-m}^{(l)}
\\& \times \binom {l-1}{m}(x+j)^m
\\&=x\sum_{l=1}^{n}\sum_{j=0}^{n} \binom {n}{l}(1-\lambda)^{-n}\binom {n}{j}
(-\lambda)^{n-j}2^l \frac{(n-1)!}{(l-1)!}\sum_{m=0}^{l-1} B_{l-1-m}^{(l)}\binom {l-1}{m}
\\&\times \left(\frac{1-\lambda}{e^t-\lambda}\right)^n (x+j)^m
\\&=\frac{x(n-1)!}{(1-\lambda)^n}\sum_{l=1}^{n}\sum_{j=0}^{n} \sum_{m=0}^{l-1}
\frac{\binom {n}{l}\binom {n}{j} \binom {l-1}{m}(-\lambda)^{n-j}2^l}{(l-1)!}
 B_{l-1-m}^{(l)}H_m^{(n)}(x+j | \lambda).
\end{align*}
Recall here that, for any $g(t) \in {\mathcal{F}}$, the Pincherle derivative is given by
\begin{equation*}\tag{25}
g^{'}(t)=g(t)x-xg(t)
\end{equation*}
as linear operators on ${\mathbb{P}}$, (see [22]).

By (18), (24) and (25), we get
\begin{align*}\tag{26}
D_n(x|\lambda)
&=\frac{e^t-\lambda}{1-\lambda}M_n(x)
\\&=\frac{(n-1)!}{(1-\lambda)^n}\sum_{l=1}^{n}\sum_{j=0}^{n} \sum_{m=0}^{l-1}\frac{\binom {n}{l}\binom {n}{j} \binom {l-1}{m}(-\lambda)^{n-j}2^l}{(l-1)!}
 B_{l-1-m}^{(l)} \\
 &\times \left\{xH_m ^{(n-1)}(x+j|\lambda)+\frac{1}{1-\lambda}H_m ^{(n)} (x+j+1|\lambda)\right\}.
\end{align*}
Therefore, by (26), we obtain the following theorem.
\par\bigskip
{\bf Theorem 2  }. For $n \geq 1$, we have
\begin{align*}
&D_n(x|\lambda)\\&=\frac{(n-1)!}{(1-\lambda)^n}\sum_{l=1}^{n}\sum_{j=0}^{n} \sum_{m=0}^{l-1}
\frac{\binom {n}{l}\binom {n}{j} \binom {l-1}{m}(-\lambda)^{n-j}2^l}{(l-1)!}
 B_{l-1-m}^{(l)}x \\
& \times \left\{xH_m ^{(n-1)}(x+j|\lambda)+\frac{1}{1-\lambda}H_m ^{(n)} (x+j+1|\lambda)\right\}.
\end{align*}
\par\bigskip
{\bf Remark  }. (a) Let us consider the following associated sequence for $\frac{t}{e^t+1}$:
\begin{align*}\tag{27}
S_n(x) \sim \left(1,\frac{t}{e^t+1}\right).
\end{align*}
By (27), we get
\begin{align*}\tag{28}
S_n(x)=x(e^t+1)^nx^{-1}x^n=x  \sum_{j=0}^{n}\binom {n}{j}(x+j)^{n-1}.
\end{align*}
From (14), (17) and (27), we have
\begin{align*}\tag{29}
M_n(x)=x\left(\frac{t}{e^t-1}\right)^nx^{-1}S_n(x)=\sum_{j=0}^{n}\binom {n}{j}xB_{n-1}^{(n)}(x+j).
\end{align*}
Thus, by (18) and (29), we get
\begin{align*}\tag{30}
 D_n(x|\lambda)
&=\left(\frac{e^t-\lambda}{1-\lambda}\right)M_n(x)
\\&=\frac{1}{1-\lambda}
\sum_{j=0}^{n}\binom {n}{j}(x+1)B_{n-1}^{(n)}(x+j+1)+\frac{\lambda}{\lambda-1}\sum_{j=0}^{n}\binom {n}{j}xB_{n-1}^{(n)}(x+j).
\end{align*}
\\(b) From (17), we note that we can obtain another expression  of $D_n(x|\lambda)$ as follow :
 \begin{align*}\tag{31}
 D_n(x|\lambda)
&=\left(\frac{e^t-\lambda}{1-\lambda}\right)M_n(x)
\\&=\sum_{k=0}^{n}\binom {n}{k}(n-1)_{n-k}2^k\frac{e^t-\lambda}{1-\lambda}(x)_k
\\&=\frac{1}{1-\lambda}
\sum_{k=0}^{n}\binom {n}{k}(n-1)_{n-k}2^k \left\{(x+1)_k - \lambda(x)_k\right\}.
\end{align*}

For $\lambda \in \mathbb{C}$ with $\lambda \neq 1$, Changhee polynomials of order $a$ are
defined by the generating function to be
\begin{align*}\tag{32}
\sum_{k=0}^{\infty}C_k^{(a)}(x|\lambda)\frac{t^k}{k!}=\left(\frac{t+1-\lambda}{1-\lambda}\right)^a(1+t)^{x}, \quad({see \,[15,23,24,25,26]}).
\end{align*}
In the specials case, $x=0$, $C_k^{(a)}(0| \lambda)=C_k^{(a)}(\lambda)$ are called $n$-th Changhee numbers of order $a$.
From (12) and (21), we note that
\begin{align*}\tag{33}
C_n^{(a)}(x|\lambda) \sim \left(\left(\frac{1-\lambda}{e^t-\lambda}\right)^a,e^t-1\right), \quad a\neq0 .
\end{align*}
Thus, by (33), we get
\begin{align*}\tag{34}
\left(\frac{e^t-\lambda}{1-\lambda}\right)^aC_n^{(a)}(x|\lambda) \sim \left(1,e^t-1\right) .
\end{align*}
From (14), $x^n \sim (1,t)$ and (35), we have
\begin{align*}\tag{35}
\left(\frac{e^t-\lambda}{1-\lambda}\right)^aC_n^{(a)}(x|\lambda)&=x\left(\frac{t}{e^t-1}\right)^nx^{-1}x^n
\\&=x\left(\frac{t}{e^t-1}\right)^nx^{n-1}=xB_{n-1}^{(n)}(x).
\end{align*}
Thus, by (14), (22), (25) and (35), we get
\begin{align*}\tag{36}
C_n^{(a)}(x|\lambda)&=\left\{x\left(\frac{1-\lambda}{e^t-\lambda}\right)^a-\frac{a}{1-\lambda}e^t\left(\frac{1-\lambda}{e^t-\lambda}\right)^{a+1}\right\}B_{n-1} ^{(n)}(x)
\\&=\sum_{l=0} ^{n-1}\binom {n-1}{l}B_l ^{(n)}\left\{x\left(\frac{1-\lambda}{e^t-\lambda}\right)^a-\frac{a}{1-\lambda}e^t\left(\frac{1-\lambda}{e^t-\lambda}\right)^{a+1}\right\}x^{n-1-l}
\\&=\sum_{l=0}^{n-1}\binom {n-1}{l}B_l^{(n)}\left\{xH_{n-1-l}^{(a)}(x|\lambda) -\frac{a}{1-\lambda}H_{n-1-l} ^{(a+1)} (x+1|\lambda)\right\}.
\end{align*}
Therefore, by (36), we obtain the following theorem.
\par\bigskip
{\bf Theorem 3  }. For $n \in \mathbb{N}$, we have
\begin{align*}
C_n^{(a)}(x|\lambda)=\sum_{l=0}^{n-1}\binom {n-1}{l}B_l^{(n)}\left\{xH_{n-1-l}^{(a)}(x|\lambda) -\frac{a}{1-\lambda}H_{n-1-l} ^{(a+1)} (x+1|\lambda)\right\}.
\end{align*}

\vskip .1in

Let
\begin{align*}\tag{37}
p(x) \sim \left(1,t\left(\frac{e^t-\lambda}{1-\lambda}\right)^a\right).
\end{align*}
By $x^n \sim (1,t)$ and (14), we get
\begin{align*}\tag{38}
p(x)&=x\left(\frac{t}{t\left(\frac{e^t-\lambda}{1-\lambda}\right)^a}\right)^n x^{-1} x^n
\\&=x\left(\frac{1-\lambda}{e^t-\lambda}\right)^{an}x^{n-1}=xH_{n-1}^{(an)}(x |\lambda),
\end{align*}
where $a \in \mathbb{Z}_+$.
\\Therefore, by (38), we obtain the following theorem.
\par\bigskip
{\bf Theorem 4 }. For $n \in \mathbb{N}$ and $a \in \mathbb{Z}_+$, we have
\begin{align*}
xH_{n-1}^{(an)}(x|\lambda) \sim \left(1,t\left(\frac{1-\lambda}{e^t-\lambda}\right)^a\right).
\end{align*}
\par\bigskip
By the same method of (38), we easily see that
\begin{align*}\tag{39}
xB_{n-1}^{(bn)}(x) \sim \left(1,t\left(\frac{e^t-1}{t}\right)^b\right), \quad b\in \mathbb{Z}_+.
\end{align*}
For $n\geq 1$ , by (14), Theorem 4 and (39), we get
\begin{align*}\tag{40}
xH_{n-1}^{(an)}(x|\lambda)&=x\left(\frac{t\left(\frac{e^t-1}{t}\right)^b}{t\left(\frac{e^t-\lambda}{1-\lambda}\right)^a}\right)^n
x^{-1}xB_{n-1}^{(an)}(x)
 \\&=x\left(\frac{1-\lambda}{e^t-\lambda}\right)^{an}\left(\frac{e^t-1}{t}\right)^{bn}B_{n-1}^{(an)}(x).
\end{align*}
Thus, from (40), we have
\begin{align*}\tag{41}
\frac{(e^t-\lambda)^{an}}{(1-\lambda)^{an}}H_{n-1}^{(an)}(x|\lambda)=\left(\frac{e^t-1}{t}\right)^{bn}B_{n-1}^{(bn)}(x).
\end{align*}
\begin{align*}\tag{42}
LHS \quad of \quad (40)&=\frac{1}{(1-\lambda)^{an}}\sum_{j=0}^{an}\binom {an}{j}(-\lambda)^{an-j}e^{jt}H_{n-1}^{(an)}(x|\lambda)
\\&=\frac{1}{(1-\lambda)^{an}}\sum_{j=0}^{an}\binom {an}{j}(-\lambda)^{an-j}H_{n-1}^{(an)}(x+j|\lambda),
\end{align*}
and
\begin{align*}\tag{43}
RHS \quad of \quad (41)&=\left(\frac{e^t-1}{t}\right)^{bn}B_{n-1}^{(bn)}(x)
\\&=\frac{1}{t^{bn}}(bn)!\sum_{j=bn}^{\infty}S_2(j,bn)\frac{t^j}{j!}B_{n-1}^{(bn)}(x)
\\&=(bn)!\sum_{j=0}^{n-1}\frac{S_2(j+bn,bn)}{(j+bn)!}t^jB_{n-1}^{(bn)}(x)
\\&=\sum_{j=0}^{n-1}\frac{(bn)!}{(j+bn)!}\frac{(n-1)!}{(n-j-1)!}S_2(j+bn,bn)B_{n-1-j}^{(bn)}(x)
\\&=(n-1)!\sum_{j=0}^{n-1}\frac{(bn)!}{(j+bn)!}\frac{S_2(j+bn,bn)}{(n-j-1)!}B_{n-1-j}^{(bn)}(x).
\end{align*}
Therefore, by (41), (42) and (43), we obtain the following theorem.
\par\bigskip
{\bf Theorem 5 }. For $a,b \in \mathbb{Z}_+$ and $n \in \mathbb{N}$, we have
\begin{align*}
&\sum_{j=0}^{an}\binom {an}{j}(-\lambda)^{an-j}H_{n-1}^{(an)}(x+j|\lambda)
\\&=(1-\lambda)^{an}(n-1)!\sum_{j=0}^{n-1}\frac{(bn)!}{(j+bn)!(n-j-1)!}S_2(j+bn,bn)B_{n-1-j}^{(bn)}(x).
\end{align*}
\par\bigskip
{\bf Remark }.
Let us consider the following Sheffer sequence:
\begin{align*}\tag{44}
S_{n}(x|\lambda)=\left( \frac{1-\lambda}{e^t-\lambda},t\left(\frac{1-\lambda}{e^t-\lambda}\right)\right).
\end{align*}
Thus, by (44), we get
\begin{align*}\tag{45}
\frac{1-\lambda}{e^t-\lambda}S_{n}(x|\lambda) \sim \left(1,t\left(\frac{1-\lambda}{e^t-\lambda}\right)\right).
\end{align*}
By Theorem 5 and (45), we get
\begin{align*}
S_{n}(x|\lambda)&=\left( \frac{e^t-\lambda}{1-\lambda}\right)xH_{n-1}^{(n)}(x|\lambda)
\\&=\frac{1}{1-\lambda}\{ (x+1)H_{n-1}^{(n)}(x+1|\lambda)-\lambda x H_{n-1}^{(n)}(x|\lambda) \}.
\end{align*}
The Daehee polynomials of the second kind are defined by the generating function to be
\begin{align*}\tag{46}
\sum_{k=0}^{\infty}D_k^*(x|\lambda)\frac{t^k}{k!}=\frac{1}{1-t}\left(\frac{1-\lambda t}{1-t}\right)^x,
\end{align*}
where $\lambda \in \mathbb{C}$ with $\lambda \neq 1$, $(see[2,14,26,28])$.
\\From (12) and (46), we note that
\begin{align*}\tag{47}
D_n^*(x|\lambda) \sim \left(\frac{1-\lambda}{e^t-\lambda},  \frac{e^t-1}{e^t-\lambda}\right).
\end{align*}
Let us consider the $\lambda$-analogues of Mittag-Leffler sequence as follows:
\begin{align*}\tag{48}
M_n^*(x|\lambda) \sim \left(1,  \frac{e^t-1}{e^t-\lambda}\right).
\end{align*}
From (14), $x^n \sim(1,t)$ and (48), we have
\begin{align*}\tag{49}
M_n^*(x|\lambda)&=x\left(\frac{e^t-\lambda}{e^t-1}t\right)^n x^{n-1}
=x(e^t-\lambda)^n\left(\frac{t}{e^t-1}\right)^{n}x^{n-1}
\\&=x\sum_{l=0}^{n}\binom {n}{l}(-\lambda)^{n-l}e^{lt}B_{n-1}^{(n)}(x)
\\&=x\sum_{l=0}^{n}\binom {n}{l}(-\lambda)^{n-l}B_{n-1}^{(n)}(x+l).
\end{align*}

Therefore, by (49), we obtain the following theorem.
\par\bigskip
{\bf Theorem 6 }. For $n \in \mathbb{N}$, let $M_n^*(x|\lambda) \sim \left(1,\frac{e^t-1}{e^t-\lambda}\right)$
\\Then we have
\begin{align*}
M_n^*(x|\lambda)=\sum_{l=0}^{n}\binom {n}{l}(-\lambda)^{n-l}xB_{n-1}^{(n)}(x+l).
\end{align*}
\par\bigskip
From (13) and Theorem 6, we have
\begin{align*}\tag{50}
M_n^*(x|\lambda)=x\sum_{l=0}^{n}\sum_{j=0}^{n-1}\binom {n}{l}\binom {n-1}{j}(-\lambda)^{n-l}B_{n-1-j}^{(n)}(x+l)^j.
\end{align*}
By (47) and (48), we get
\begin{align*}\tag{51}
D_n^*(x|\lambda)&=\left(\frac{e^t-\lambda}{1-\lambda}\right) M_n^*(x|\lambda)
\\&=\frac{1}{1-\lambda }\sum_{j=0}^{n}\binom {n}{j}\{(x+1)B_{n-j}^{(n)}(x+1+j)-\lambda xB_{n-1}^{(n)}(x+j)\}(-\lambda)^{n-j}.
\end{align*}
Let us consider the following associated sequences:
\begin{align*}\tag{52}
T_n^*(x) \sim \left(1, \frac{2t}{1+t^2}\right).
\end{align*}
From (12) and (52), we have
\begin{align*}\tag{53}
\sum_{k=0}^{\infty}T_k^*(x)\frac{t^k}{k!}=exp\left(x\left(\frac{1-\sqrt{1-t^2}}{t}\right)\right).
\end{align*}
By (14), (52) and $x^n \sim (1,t)$, we get
\begin{align*}
T_n^*(x)&=x\left(\frac{t}{\frac{2t}{1+t^2}}\right)^nx^{-1}x^n
=x\left(\frac{1+t^2}{2}\right)^nx^{n-1}
\\&=\left(\frac{1}{2}\right)^nx\sum_{l=0}^{[\frac{n-1}{2}]}\binom {n}{l}t^{2l}x^{n-1}
\\&=\left(\frac{1}{2}\right)^n\sum_{l=0}^{[\frac{n-1}{2}]}\binom {n}{l} \frac{(n-1)!}{(n-2l-1)!}x^{n-2l}.
\end{align*}

\par\bigskip

\par\noindent
\mpn { \bpn {\small Dae San {\sc KIM} \mpn
 Department of Mathematics,\pn
Sogang University, Seoul 121-742, S.Korea
  \pn {\it E-mail:}\ {\sf dskim@sogang.ac.kr} }

\mpn { \bpn {\small Taekyun {\sc KIM} \mpn
 Department of Mathematics,\pn
Kwangwoon University, Seoul 139-701, S.Korea
  \pn {\it E-mail:}\ {\sf tkkim@kw.ac.kr} }

  \mpn { \bpn {\small Seog-Hoon {\sc Rim} \mpn Department of Mathematics
Education, \pn Kyungpook National University, Daegu 702-701, S.
Korea \pn {\it E-mail:}\ {\sf shrim@knu.ac.kr} }
 \


\begin{thebibliography}{9}


\bibitem{[1]} S. Araci, M. Acikgoz, {\it A note on the Frobenius-Euler numbers and polynomials
associated with Bernstein polynomials}, Adv. Stud. Contemp. Math. 22 (2012), no. 3, 399-406.
\bibitem{[1]} S. Araci, M. Acikgoz, A. Esi, {\it A note on the q-Dedekind-type Daehee-Changhee sums with weight alpha arising from modified q-Genocchi polynomials with weight alpha}, arXiv:1211.2350
\bibitem{[2]} S. Araci, M. Acikgoz, {\it Extended q-Dedekind-type Daehee-Changhee sums associated with Extended q-Euler polynomials},  arXiv:1211.1233
\bibitem{[3]} S. Araci, E Sen, M. Acikgoz, {\it A note on the modified q-Dedekind sums},  arXiv:1212.5837
\bibitem{[4]} L. Carlitz, {\it Some generating functions for Laguerre polynomials}, Duke Math. J. 35  (1968) 825-827.
\bibitem{[5]} J. Choi, {\it A note on $p$-adic integrals associated with Bernstein and $q$-Bernstein polynomials}, Adv. Stud. Contemp. Math. 21 (2011), no. 2, 133-138.
 \bibitem{[6]} R.Dere, Y.Simsek, {\it Applications of umbral algebra to some special polynomials}, Adv. Stud. Contemp. Math.  22  (2012),  no. 3, 433-438.
\bibitem{[7]} Q. Fang, T. Wang, {\it Umbral calculus and invariant sequences}, Ars Combin. 101  (2011), 257-264.
\bibitem{[8]} D. S. Kim, T.Kim, S.-H. Lee, S.-H. Rim, {\it Frobenius-Euler polynomials and umbral calculus in the p-adic case}, Advances in Difference Equations 2012, 2012:222.
\bibitem{[9]} D. S. Kim, T. Kim, {\it Some new identities of Frobenius-Euler numbers and polynomials},  Journal of Inequalities and Applications 2012, 2012:307
\bibitem{[10]} D. S. Kim, T. Kim, {\it Applications of Umbral Calculus Associated with $p$-Adic Invariant Integrals on $\Bbb Z_p$}, Abstract and Applied Analysis 2012 (2012), Article ID 865721, 12 pages
\bibitem{[11]} D. S. Kim, T. Kim, {\it Some identities of Frobenius-Euler polynomials arising from umbral calculus}, Advances in Difference Equations 2012, 2012:196 .
\bibitem{[12]} T. Kim, S.-H. Rim, D. V. Dolgy, S.-H. Lee, {\it Some identities on Bernoulli and Euler polynomials arising from the orthogonality of Laguerre polynomials},  Advances in Difference Equations 2012, 2012:201.
\bibitem{[13]} T. Kim, {\it Some identities on the $q$-Euler polynomials of higher order and $q$-Stirling numbers by the fermionic $p$-adic integral on $\Bbb Z_p$}, Russ. J. Math. Phys. 16  (2009),  no. 4, 484-491.
\bibitem{[14]} T. Kim, {\it An invariant p-adic integral associated with Daehee numbers}, Integral Transforms Spec. Funct. 13 (2002), no. 1, 65-69.
\bibitem{[15]} T. Kim, {\it $p$-adic $q$-integrals associated with the Changhee-Barnes' $q$-Bernoulli polynomials}, Integral Transforms Spec. Funct. 15 (2004), no. 5, 415-420.
\bibitem{[16]} T. Kim, Y. Simsek, {\it Analytic continuation of the multiple Daehee q-l-functions associated with Daehee numbers}, Russ. J. Math. Phys. 15  (2008),  no. 1, 58-65
\bibitem{[17]} H. Ozden, I. N. Cangul, Y. Simsek, {\it Remarks on $q$-Bernoulli numbers associated with Daehee numbers}, Adv. Stud. Contemp. Math.  18  (2009),  no. 1, 41-48.
\bibitem{[18]} Y. Simsek,  S.-H. Rim, L.-C. Jang, D.-J. Kang, J.-J. Seo, {\it A note on q-Daehee sums}, J. Anal. Comput. 1  (2005), no. 2, 151-160.
\bibitem{[19]} T. J. Robinson, {\it Formal calculus and umbral calculus}, Electron. J. Combin. 17  (2010), no. 1, Research Paper 95, 31 pp.
\bibitem{[20]} S. Roman, {\it More on the umbral calculus, with emphasis on the q-umbral calculus}, J. Math. Anal. Appl.107 (1985), 222-254. MR0786026 (86h:05024)
\bibitem{[21]} S. Roman, {\it The umbral calculus}, Dover Publ. Inc. New York, 2005.
\bibitem{[22]} C. Ryoo, {\it Some relations between twisted q-Euler numbers and Bernstein polynomials}, Adv. Stud. Contemp. Math.  21  (2011),  no. 2, 217-223.
\bibitem{[23]} C. S. Ryoo, D.-W. Park, S,-H. Rim, {\it On the real roots of the Changhee-Barnes' $q$-Bernoulli polynomials}, JP J. Algebra Number Theory Appl. 5 (2005),  no. 2, 293-305.
\bibitem{[24]} C. S. Ryoo, T. Kim, R. P. Agarwal, {\it Exploring the multiple Changheeq-Bernoulli polynomials}, Int. J. Comput. Math. 82  (2005), no. 4, 483-493.
\bibitem{[25]} C. S. Ryoo, H. Song, {\it On the real roots of the Changhee-Barnes' q-Bernoulli polynomials}, Proceedings of the 15th International Conference of the Jangjeon Mathematical Society,  63-85, Jangjeon Math. Soc., Hapcheon, 2004.
\bibitem{[26]} Y. Simsek, {\it Special functions related to Dedekind-type DC-sums and their applications}, Russ. J. Math.  Phys. 17 (2010) 495-508.
\bibitem{[27]} Y. Simsek, {\it Generating functions of the twisted Bernoulli numbers and polynomials associated with their interpolation functions}, Adv. Stud. Contemp. Math.  16  (2008),  no. 2, 251-278.
\bibitem{[28]} Y. Simsek,  I.-S. Pyung, {\it Barnes' type multiple Changheeq-zeta functions}, Adv. Stud. Contemp. Math.  10  (2005),  no. 2, 121-129.

\end{thebibliography}
\end{document}